\documentclass[preprint]{elsarticle}
\usepackage{tikz-cd}
\usepackage{cmap} % for russuan cope paste
           \usepackage[T2A]{fontenc}   % Or any other encoding.
           \usepackage[utf8]{inputenc}
           \usepackage[english]{babel}
\usepackage{amsmath}
\usepackage{amsfonts}
\usepackage{amssymb}
\usepackage{amsthm}
\usepackage[pdftex,unicode]{hyperref}
\usepackage{indentfirst} % включить отступ у первого абзаца
\usepackage{amscd}
\usepackage[all,cmtip]{xy}
\usepackage{mathtools}
\usepackage{comment}
\usepackage{tikz-cd}
\usetikzlibrary{babel}

  \def\R{\mathbb R}

\def\om{\omega}

\def\dim{\operatorname{dim}}
\def\al{\alpha}
\def\ga{\gamma}

\def\ph{\varphi}

\def\nfs/{NFS}
\def\cdp/{CDP}
\def\cdpz/{CDP${}_0$}

\def\es{\varnothing}

\def\sset#1{\{#1\}}

\def\set#1{\bbset#1\eeset}
\def\bbset#1:#2\eeset{\{#1\,:\,#2\}}

\def\bbsett#1:#2\eesett{\{#1\,:\,\text{#2}\}}

\def\ibbset#1:#2\ieeset{(#1)_{#2}}

\def\DD#1#2#3{\D_{#1,#2}(#3)}

\def\cF{{\mathcal F}}

\def\cU{{\mathcal U}}

\def\cV{{\mathcal V}}

\newcommand\restrA[2]{{% we make the whole thing an ordinary symbol
  \left.\kern-\nulldelimiterspace % automatically resize the bar with \right
  #1 % the function
  \vphantom{\big|} % pretend it's a little taller at normal size
  \right|_{#2} % this is the delimiter
  }}

\newcommand\restrB[2]{\ensuremath{\left.#1\right|_{#2}}}

\def\restr#1#2{\restrB{#1}{#2}}

\def\pwr#1_#2{#1^{[#2]}}

\def\term#1{{\it #1}}

\def\pp{{\mathsf{P}}}

\def\dddm#1(#2){N_{#1}(#2)}
\def\dddb#1(#2){B_{#1}(#2)}

\def\et(#1){ (#1)}

\def\bitem#1,#2.{ $#2\nrightarrow #1$:\ }
\newtheorem{assertion}{Statement}
\newtheorem{proposition}{Proposition}
\newtheorem{theorem}{Theorem}
\newtheorem{lemma}{Lemma}
\newtheorem*{lemma*}{Lemma}
\newtheorem{cor}{Corollary}

\theoremstyle{definition}
\newtheorem{definition}{Definition}

\def\oo#1/{$O_{#1}$}

\def\gd/{$G_\delta$}

\def\pp{\Psi}

\def\DD{\operatornamewithlimits{%
  \mathchoice{\vcenter{\hbox{\huge \Delta}}}
             {\vcenter{\hbox{\Large \Delta}}}
             {\mathrm{\Delta}}
             {\mathrm{\Delta}}}}

\def\DD{\operatornamewithlimits{\Delta}}

\def\DD{\operatornamewithlimits{\mbox{\hbox{\huge \Delta}}}}

\def\DD{\operatornamewithlimits{%
  \mathchoice{\vcenter{\hbox{\huge $\Delta$}}}
             {\vcenter{\hbox{\Large $\Delta$}}}
             {\mathrm{\Delta}}
             {\mathrm{\Delta}}}}

\def\si{\sigma}

\def\ssi#1#2{#1_{[#2]}}
\def\ssr#1#2{#1_{[#2)}}
\def\ssl#1#2{#1_{(#2]}}
\def\ppn#1{$(\star_{#1})$}
\def\pp/{\ppn 1}
\def\pppp/{\ppn 2}
\def\bb{\mathcal{B}}
\def\bz{B}

\def\sL{\mathsf{L}}
\def\sR{\mathsf{R}}
\def\sg#1#2{\mathsf{{#1}_{#2}}}

\def\cls#1,#2/{\mathsf{#1}_{#2}}

\begin{document}

\begin{frontmatter}

\title{The product of Lindel\"of groups and $\R$-factorizability}
\author{Evgeny Reznichenko}

\ead{erezn@inbox.ru}

\address{Department of General Topology and Geometry, Mechanics and Mathematics Faculty, M.~V.~Lomonosov Moscow State University, Leninskie Gory 1, Moscow, 199991 Russia}

\begin{abstract}
Lindel\"of topological groups $\sg G1$, $\sg H1$, $\sg G2$, $\sg H2$ are constructed in such a way that the products of $\sg G1\times \sg H1$ and $\sg G2\times \sg H2$ are not $\R$-factorizable groups and (1) the group $\sg G1\times \sg H1$ is not pseudo-$\aleph_1$-compact; (2) the group $\sg G2\times \sg H2$ 
is a separable not  normal group and
 contains a discrete closed subset of the cardinality $2^\om$.
\end{abstract}
\end{frontmatter}

\section{Introduction}
\label{sec:intro}
A topological group $G$ is said to be $\R$-factorizable if for any continuous function $f: G\to \R$ there exists a topological separable metrizable group $G$, a continuous homomorphism $\ph: G\to H$, and a continuous function $h: H\to \R$, such that $f=h\circ \ph$, i.e. the following diagram is commutative
\[
\begin{tikzcd}
G \ar[rr,"f"]\ar[rd,"\ph"]&& H
\\
& H \ar[ru,"h"]
\end{tikzcd}
\]
$\R$-factorizable groups play a big role in topological group theory, see Chapter 8 in \cite{at2009}.
One of the main questions in the theory of $\R$-factorizable groups:
is the product of two $\R$-factorizable groups a $\R$-factorizable group \cite{Tkachenko1991,Tkachenko2004}? Lindelöf topological groups are $\R$-factorizable \cite[Theorem 8.1.6.]{at2009}. In this note, Lindel\"of groups are constructed whose product is not $\R$-factorizable.
These examples also answer Problems 8.5.3, 8.5.4 and 8.5.9 from \cite{at2009}.

Lindel\"of topological groups $\sg G1$, $\sg H1$, $\sg G2$, $\sg H2$ are constructed in such a way that the products of $\sg G1\times \sg H1$ and $\sg G2\times \sg H2$ are not $\R$-factorizable groups and (1) the group $\sg G1\times \sg H1$ is not pseudo-$\aleph_1$-compact; (2) the group $\sg G2\times \sg H2$ 
is a separable not  normal group and
 contains a discrete closed subset of the cardinality $2^\om$ (Theorem \ref{t:products:1}).
The group $\sg G2\times \sg H2$ is a separable non-$\R$-factorizable group.
The first such group was constructed in  \cite{ReznichenkoSipacheva2013}.

Recently O.V. Sipacheva \cite{Sipacheva2023,Sipacheva2023-2} also independently found a non-$\R$-factorizable product of Lindelöf groups. Moreover, one of these groups is metrizable separable. In Section \ref{sec:dim}, products of Lindel\"of groups from \cite{Sipacheva2023,Sipacheva2023-2} are studied and it is proved that the product of Lindel\"of groups from \cite{Sipacheva2022} is also not $\R$-factorizable.

\section{Notation and definition}
\label{sec:defs}

Let $X$ be a set and $\tau_1$, $\tau_2$, ..., $\tau_n$ be a family of subsets of $X$.
Denote by $\bb(\tau_1, \tau_2, ..., \tau_n)$ the topology on $X$ whose prebase is formed by the family $\bigcup_{i=1}^n \tau_i$.
If $\tau_i$ is a topology on $X$ for $i=1,2,...,n$, then
the topology $\bb(\tau_1, \tau_2, ..., \tau_n)$ is the topology of the diagonal $\set{(x,x,...,x): x\in X}$ in the product
$
\prod_{i=1}^n (X,\tau_i).
$

Denote by $\bz(X)$ the free Boolean topological group over $X$.

A family $\cF$ is called \term{non-Archimedean} if for any $A,B\in\cF$ one of the following conditions is true: $A\cap B=\es$; $A\subset B$; $B\subset A$. 
Spaces that have a non-Archimedean base are called \term{non-Archimedean spaces}. Separable metrizable zero-dimensional spaces are non-Archimedean.

We call a space $X$ \term{co-non-Archimedean} if there is a non-Archimedean topology $\tau$ on $X$ and a family $\cF$ of  $\tau$-closed subsets such that the topology $X$ coincides with $\bb (\tau,\cF)$.

\begin{theorem}[\cite{grs1997}]\label{t:pog-defs:grs}
If $X$ is a co-non-Archimedean space, then $X$ is a retract $\bz(X)$.
\end{theorem}

The group with the topology is called
\begin{enumerate}
\item
\term{paratopological} if the multiplication in the group is (jointly) continuous;
\item
\term{semitopological} if the multiplication in the group is separately continuous.
\end{enumerate}

Consider the following classes of spaces and groups with topology:
\begin{enumerate}
\item[($\cls TS,sm/$)]
separable metrizable spaces;
\item[($\cls TS,cnw/$)]
spaces with a countable network;
\item[($\cls TG,sm/$)]
separable metrizable topological groups;
\item[($\cls PG,sm/$)]
separable metrizable paratopological groups;
\item[($\cls SG,cnw/$)]
semitopological groups with a countable network.
\end{enumerate}
Let $\cls G,/$ be some class of groups with topology.
A group with the topology $G$ is called \term{$\cls G,/$-factorizable} if for any continuous function $f: G\to \R$ there exists a group $H\in \cls G,/$, a continuous homomorphism $\ph: G\to H$ and a continuous function $h: H\to \R$, such that $f=h\circ \ph$.

A topological group $G$ is called \term{$\R$-factorizable} if $G$ is $\cls TG,sm/$-factorizable.
A paratopological group $G$ is called \term{$\R$-factorizable} if $G$ is $\cls PG,sm/$-factorizable.

Let $\cls S,/$ be some class of spaces.
We call a function on the product of spaces $f: X_1\times X_2 \times ... \times X_n \to \R$ \term{(separately) $\cls S,/$-factorizable} if there exist spaces $Y_1$, $ Y_2$, ..., $Y_n$ from $\cls S,/$, continuous mappings $g_i: X_i\to Y_i$ for $i=1,2,...,n$ and a (separately)  continuous function $ h: Y= Y_1\times Y_2 \times ... \times Y_n\to \R$, so $f=h\circ \prod_{i=1}^n g_i$, that is, the following diagram is commutative.
\[
\begin{tikzcd}
X_1\times X_2 \times ... \times X_n
\ar[rr,"f"]
\ar[dd,"\prod\limits_{i=1}^n g_i"']
&& \R
\\
\\
Y_1\times Y_2 \times ... \times Y_n
\ar[rruu,"h"']
\end{tikzcd}
\]
We call a function $f$  \term{(separately) $\R$-factorizable} if $f$ is (separately) $\cls TS,sm/$-factorizable.

A product of spaces $X=X_1\times X_2 \times ... \times X_n$ is  called \term{(separately) $\cls S,/$-factorizable} if any continuous function on $f:X\to\R$ is (separately) $\cls S,/$-factorizable.

A product of spaces $X=X_1\times X_2 \times ... \times X_n$ is  called \term{(separately) $\R$-factorizable} if $X$ is (separately) $\cls TS,sm/$- factorizable.

\term{The discrete Suslin number} $dc(X)$ of the space $X$ is called the supremum of cardinals $\tau$,
such that there exists a discrete family of open sets of cardinality $\tau$.
Spaces with a countable discrete Suslin number are also called \term{pseudo-$\aleph_1$-compact}.

%\term{Extend} $e(X)$ of $X$ is called the supremum of cardinals $\tau$, such that there exists a discrete closed set of cardinality $\tau$.

\begin{theorem}[\cite{Jones1937}]\label{t:pog-defs:Jones1937}
If $X$ is a normal separable space then 
the space does not contain a discrete closed subset of the cardinality of the continuum.
\end{theorem}

There are two definitions of covering dimension, in the sense of \v Cech and in the sense of
Kat\v etov; following \cite{Charalambous2019}, we denote the Kat\v etov dimension by $\dim_0$. Recall
that, given a topological space $X$, $\dim_0 X$ is the least integer $n > -1$ such
that any finite cozero cover of $X$ has a finite cozero refinement of
order $n$, provided that such an integer exists. If it does not exist, then
$\dim_0 X = \infty$. A space $X$ for which $\dim_0 X = 0$ is said to be \term{strongly zero-dimensional}.

In \cite{EngelkingBookGT}, as in many other places, Kat\v etov dimension is denoted as $\dim$.

\begin{theorem}\label{t:pog-defs:rdim}
Let the topological group $H$ be $\R$-factorizable.
\begin{enumerate}
\item
If $H$ is embedded in a topological group $G$, then $\dim_0 H\leq \dim_0 G$ \cite[Theorem 2.7]{Tkachenko1991}.
\item
If $H$ is zero-dimensional then $H$ is strongly zero-dimensional \cite[Theorem 3.3]{Shakhmatov1990}.
\end{enumerate}
\end{theorem}

\begin{definition}[Definition 2.28 \cite{KozlovPasynkov2007}]
A product space $X \times Y$ is said to be \term{rectangular} if any ﬁnite cozero-cover $\cU$ of the space $X\times Y$ has a $\si$-locally ﬁnite reﬁnement $\cV $ consisting of sets of the form $U\times V$, where $U$ is cozero subset of $X$ and $V$ is cozero subset of $Y$.
\end{definition}

\begin{theorem}[Theorem 2.32 \cite{KozlovPasynkov2007}]\label{t:pog-defs:rect}
If a product space $X \times Y$ is rectangular, then
\[
\dim_0 X\times Y \leq \dim_0 X + \dim_0 Y.
\]
\end{theorem}

\section{Product of groups}
\label{sec:products}
Let $X\subset \R$ and $A\subset X$, $\tau$ be the topology of $X$ induced from $\R$.
Denote
\begin{align*}
\ssi\tau A&=\bb(\tau,\set{\sset x:x\in X\setminus A}),
&
\ssi XA&=(X,\ssi\tau A),
\\
\ssr\tau A&=\bb(\tau,\set{X\cap [x,+\infty):x\in X\setminus A}),
&
\ssr XA&=(X,\ssr\tau A),
\\
\ssl\tau A&=\bb(\tau,\set{X\cap (-\infty,x]:x\in X\setminus A}),
&
\ssl XA&=(X,\ssl\tau A).
\end{align*}
The following assertions follow from the definition.
\begin{assertion}\label{a:products:1}
The spaces $\ssl XA$ and $\ssr XA$ are continuous images of $\ssi XA$.
\end{assertion}
\begin{assertion}\label{a:products:2}
If $X\setminus A$ is dense in $\R$, then the spaces $\ssi XA$, $\ssl XA$, $\ssr XA$ are co-non-Archimedean.
\end{assertion}

There exists a partition $\set{M_\al:\al<2^\om}$ of the line $\R$, such that $\left(\ssi \R{M_\al}\right)^\om$ is Lindel\"of for $ \al<2^\om$ \cite{OkunevTamano1996}.

Let $A=M_0$, $B=M_1$, $C=\R\setminus(A\cup B)$, $\sL= A\cup C$, $\sR=B\cup C$.

The construction and Statements \ref{a:products:1} and \ref{a:products:2} imply the following statement.
\begin{proposition}\label{p:products:1}
The spaces $\ssi{\sL}{A}$, $\ssl{\sL}{A}$, $\ssr{\sL}{A}$,
$\ssi{\sR}{B}$, $\ssl{\sR}{B}$, $\ssr{\sR}{B}$ are co-non-Archimedean and  Lindel\"of to countable degree.
\end{proposition}

\begin{lemma}\label{l:products:1}
If a separately continuous function $f: X\times Y\to \R$ is separately $\cls TS,cnw/$-factorizable, then $\ph(X)$ has a countable network, where
\[
\ph: X\to C_p(Y),\ \ph(x)(y)=f(x,y).
\]
\end{lemma}
\begin{proof}
Since $f$ is a separately continuous  $\cls TS,cnw/$-factorizable function, there exist
spaces $X_*$ and $Y_*$ with a countable network,
continuous mappings $g_X: X\to X_*$, $g_Y: Y\to Y_*$ and a separately continuous function $h: X_*\times Y_*\to \R$, such that
$f= h\circ (g_X \times g_Y)$. Then
\[
\psi: X_* \to C_p(Y),\ \psi(x_*,y) = h(x_*,g_Y(y))
\]
is continuous and $\ph=g_X\circ\psi$.
Then $\ph(X)$ is a continuous image of $X_*$ and hence $\ph(X)$ has a countable network.
\end{proof}

\begin{proposition}\label{p:products:2}
%{\rm(1)} The spaces $\ssl{\sL}{A}$ and $\ssr{\sR}{B}$ are separable.
{\rm(1)}
The product $\ssi{\sL}{A} \times \ssi{\sR}{B}$ is not Lindel\"of, has a discrete Suslin number $2^\om$, and is not separately $\cls TS,cnw/$-factorizable.
{\rm(2)}
The product $\ssl{\sL}{A} \times \ssr{\sR}{B}$ is separable, contains a discrete closed subset of the cardinality $2^\om$, is not normal, and is not separately $\cls TS,cnw/$-factorizable.
\end{proposition}
\begin{proof}
Denote $D=\set{(x,x): (x,x)\in \sL\times \sR}=\set{(x,x): x\in C}$.

(1) The family $\set{\sset p: p\in D}$ is a discrete family of open sets in $\ssi{\sL}{A} \times \ssi{\sR}{B}$.
Let us define the function $f: \sL\times \sR\to \R$,
\[
f(x,y)=\begin{cases}
0,& \text{if }y\geq x,
\\
1,& \text{if }y < x.
\end{cases}
\]
Let us put
\[
\ph: \sL\to C_p(\sR),\ \ph(x)(y)=f(x,y).
\]
Then
\[
\ph(x)(y)=\begin{cases}
0,& \text{if }y\geq x,
\\
1,& \text{if }y < x.
\end{cases}
\]
The uncountable space $\ph(C)$ embeds in the Songefrey line and, therefore, $\ph(C)$ and $\ph(\sL)$ do not have a countable network.

The function $f$ is continuous with respect to the topology $\ssi{\sL}{A} \times \ssi{\sR}{B}$.
Since $\ph(\sL)$ does not have a countable network, it follows from the Lemma \ref{l:products:1} that the function $f$ is not separately $\cls TS,cnw/$-factorizable.

(2)
The topologies $\ssl{\sL}{A}$ and $\ssr{\sR}{B}$ are weaker than the Songegfrey topology, and the Songegfrey line is a hereditarily separable space. Hence the spaces $\ssl{\sL}{A}$, $\ssr{\sR}{B}$ and $\ssl{\sL}{A}\times \ssr{\sR}{B}$ are separable.

The set $D$ is a discrete closed subset of $\ssl{\sL}{A}\times \ssr{\sR}{B}$ and $|D|=2^\om$. 
Theorem \ref{t:pog-defs:Jones1937} implies that the product $\ssl{\sL}{A}\times \ssr{\sR}{B}$ is not normal.
The function $f$ is continuous with respect to the topology $\ssl{\sL}{A} \times \ssr{\sR}{B}$.
Since $\ph(\sL)$ does not have a countable network, it follows from the Lemma \ref{l:products:1} that the function $f$ is not separately $\cls TS,cnw/$-factorizable.
\end{proof}

\begin{proposition}\label{p:products:3}
Let $\cls T,/$ be a class of spaces closed with respect to taking subspaces,
  % and final products,
$\cls PT,/$ is the class of paratopological groups in $\cls T,/$ and
$\cls ST,/$ is the class of semitopological groups in $\cls T,/$.
Let $G$ and $H$ be groups with topology.
\begin{enumerate}
\item
If the product of groups $G\times H$ is $\cls PT,/$-factorizable, then the product of spaces $G\times H$ is $\cls T,/$-factorizable.
\item
If the product of groups $G\times H$ is $\cls ST,/$-factorizable, then the product of spaces $G\times H$ is separately $\cls T,/$-factorizable.
\end{enumerate}
\end{proposition}
\begin{proof}
Let $f: G\times H\to \R$ be a continuous function.

(1) Since the group $G\times H$ is $\cls PT,/$-factorizable,
then there exists a paratopological group $S\in \cls T,/$, a continuous homomorphism $g: G\times H\to S$, and a continuous function $h: S\to\R$, so that $f=h\circ g $. Let $G_*=g(G)$, $H_*=g(H)$, $g_G=\restr g{G}$, $g_H=\restr g{H}$,
\[
h_*: G_*\times H_*\to \R,\ (x_*,y_*)\mapsto h(x_*y_*).
\]
Then the function $h$ is continuous and $f=h_*\circ(g_G\times g_H)$.

(2) Since the group $G\times H$ is $\cls ST,/$-factorizable,
then there exists a semitopological group $S\in \cls T,/$, a continuous homomorphism $g: G\times H\to S$, and a continuous function $h: S\to\R$, so that $f=h\circ g $. Let $G_*=g(G)$, $H_*=g(H)$, $g_G=\restr g{G}$, $g_H=\restr g{H}$,
\[
h_*: G_*\times H_*\to \R,\ (x_*,y_*)\mapsto h(x_*y_*).
\]
Then the function $h$ is separately continuous and $f=h_*\circ(g_G\times g_H)$.
\end{proof}

\begin{cor}\label{c:products:1}
Let $G$ and $H$ be groups  with topology.
If the product of groups $G\times H$ is $\cls SG,cnw/$-factorizable, then the product of spaces $G\times H$ is separately $\cls TS,cnw/$-factorizable.
\end{cor}

\begin{cor}\label{c:products:1+1}
Let $G$ and $H$ be groups  with topology.
Assume that one of the following conditions is true:
\begin{enumerate}
\item
the product of groups $G\times H$ is $\cls PG,sm/$-factorizable;
\item
the groups $G$ and $H$ are paratopological groups, and the group $G\times H$ is $\R$-factorizable.
\end{enumerate}
Then the product of spaces $G\times H$ is $\R$-factorizable.
\end{cor}

\begin{proposition}\label{p:products:4}
Let $\cls T,/$ be a class of spaces.
Let $X$ and $Y$ be topological spaces, $X'$ retract $X$ and $Y'$ retract $Y$.
If the product $X\times Y$ is (separately) $\cls T,/$-factorizable,
then the product $X'\times Y'$ is (separately) $\cls T,/$-factorizable.
\end{proposition}
\begin{proof}
Let $r_X: X\to X'$ and $r_Y: Y\to Y'$ be retractions, the function $f': X'\times Y'\to\R$ is continuous.
We put $f=f'\circ(r_X\times r_Y)$.
Since the product $X\times Y$ is (separately) $\cls T,/$-factorizable, there exist spaces $X_*,Y_*\in \cls T,/$, continuous mappings $g_X: X\to X_* $, $g_Y: Y\to Y_*$ and (separately) a continuous function $h: X_*\times Y_*\to\R$, such that $f=h\circ(g_X\times g_Y)$.
Let us put $g'_X=\restr{g_X}{X'}$ and $g'_Y=\restr{g_Y}{Y'}$. Then $f'=h\circ(g'_X\times g'_Y)$. Hence $f'$ is  a (separately) $\cls T,/$-factorizable function.
\end{proof}

\begin{cor}\label{c:products:2}
Let $X$ and $Y$ be topological spaces, $X'$ retract $X$ and $Y'$ retract $Y$.
If the product $X\times Y$ is separately $\cls TS,cnw/$-factorizable,
then the product $X'\times Y'$ is separately $\cls TS,cnw/$-factorizable.
\end{cor}

Define groups
\begin{align*}
\sg G1&= \bz(\ssi{\sL}{A}),
&
\sg H1&= \bz(\ssi{\sR}{B}),
\\
\sg G2&= \bz(\ssl{\sL}{A}),
&
\sg H2&= \bz(\ssr{\sR}{B}).
\end{align*}

\begin{theorem}\label{t:products:1}
The topological groups $\sg G1$, $\sg H1$, $\sg G2$, $\sg H2$ are Lindel\"of to countable degree.
\begin{enumerate}
\item
The product $\sg G1\times \sg H1$ is not a $\cls SG,cnw/$-factorizable (and then not $\R$-factorizable) group and has a discrete Suslin number $2^\om$.
\item
The product $\sg G2\times \sg H2$ is not a $\cls SG,cnw/$-factorizable (and then not $\R$-factorizable) non-normal separable group and has extend $2^\om$.
\end{enumerate}
\end{theorem}
\begin{proof}
It follows from Proposition \ref{p:products:1} and Theorem \ref{t:pog-defs:grs} that the groups $\sg G1$, $\sg H1$, $\sg G2$, $\sg H2 $ are Lindel\"of to countable degree and $\ssi{\sL}{A}$, $\ssi{\sR}{B}$, $\ssl{\sL}{A}$, $\ssr{\sR}{ B}$ is a retract of $\sg G1$, $\sg H1$, $\sg G2$, $\sg H2$, respectively.
From Proposition \ref{p:products:2}, Corollaries \ref{c:products:1} and \ref{c:products:2} it follows that the groups $\sg G1\times \sg H1$ and $\sg G2\times \sg H2$ are not $\cls SG,cnw/$-factorizable.

(1)
Since $\ssi{\sL}{A}\times \ssi{\sR}{B}$ is a retract of $\sg G1\times \sg H1$, then from Proposition \ref{p:products:2}(1) it follows that $\sg G1\times \sg H1$ has a discrete Suslin number $2^\om$.

(2)
Since $\ssl{\sL}{A}\times \ssr{\sR}{B}$ embeds in $\sg G2\times \sg H2$ and generates this group, from Proposition 
\ref{p:products:2}(2) implies that $\sg G1\times \sg H1$ is separable, has an extend $2^\om$, and is not normal.
\end{proof}

\section{$\R$-factorizability and dimension}
\label{sec:dim}
\begin{proposition}\label{p:dim:1}
If a product of spaces $X \times Y$ is $\R$-factorizable, then $X \times Y$ is rectangular.
\end{proposition}
\begin{proof}
Let $\cU=\sset{U_1,...,U_n}$ is a ﬁnite cozero cover of $X \times Y$.
For $i=1,...,n$, we fix a continuous function $f_i$ on $X \times Y$ such that $X\times Y \setminus U_i=f_i^{-1}(0)$.
Since $X \times Y$ is $\R$-factorizable, for $i=1,...,n$,
there exists separable metrizable spaces $X_i$ and $Y_i$, continuous mappings $p_i: X\to X_i$, and
$q_i: Y\to Y_i$ and a continuous function $h_i: X_i\times X_i\to \R$, such that $f_i=h_i\circ (p_i \times q_i)$.
Let us put
\begin{align*}
X_*&=\prod_{i=1}^n X_i,
&
Y_*&=\prod_{i=1}^n Y_i,
&
p &= \DD_{i=1}^n p_i,
&
q &= \DD_{i=1}^n q_i.
\end{align*}
Then $X_*$ and $Y_*$ are separable metrizable spaces and $U_i=(p\times q)^{-1}(V_i)$ for $i=1,...,n$ for some open $V_i \subset X_*\times Y_*$. There is a countable family $\ga_i$ consisting of rectangular open sets of the form $U\times V$, such that $V_i=\bigcup \ga_i$. Then the family
\[
\mu=\bigcup_{i=1}^n \set{ (p\times q)^{-1}(U\times V): U\times V \in \ga_i }
\]
is countable (and even more so $\si$-locally ﬁnite), consists of rectangular cozero sets, and is reﬁnement of $\cU$.
Therefore $X \times Y$ is rectangular.
\end{proof}

Corollary \ref{c:products:1+1}, Proposition \ref{p:dim:1}, and Theorem \ref{t:pog-defs:rect} imply the following assertions.

\begin{proposition}\label{p:dim:2}
Let $G$ and $H$ be groups with topology and let $G\times H$ be $\cls PG,sm/$-factorizable.
Then
\[
\dim_0 G\times H \leq \dim_0 G + \dim_0 H.
\]
\end{proposition}

\begin{theorem}\label{t:dim:1}
Let $G$ and $H$ be (para)topological groups and the (para)topological group $G\times H$ be $\R$-factorizable.
Then
\[
\dim_0 G\times H \leq \dim_0 G + \dim_0 H.
\]
\end{theorem}

In \cite{Sipacheva2022} the strongly zero-dimensional topological groups $\sg G3$ and $\sg H3$ are constructed, such that
the group $\sg G3$ is Lindel\"of, the group $\sg H3$ is separable metrizable\footnote{Group $\sg H3$  in the original paper with a countable network, but with a slight modification of the construction it is possible to make $\sg H3$ become metrizable separable}, and $\dim_0 \sg G3\times \sg H3>0$.
There are two arguments why the topological group $\sg G3\times \sg H3$ is not $\R$-factorizable.
First, because
\[
\dim_0 \sg G3\times \sg H3>0= \dim_0 \sg G3+ \dim_0 \sg H3,
\]
then the theorem \ref{t:dim:1} implies that $\sg G3\times \sg H3$ is not $\R$-factorizable.
Second, since the group $\sg G3\times \sg H3$ is zero-dimensional and is not strongly zero-dimensional, Theorem \ref{t:pog-defs:rdim}(2) implies that $\sg G3\times \sg H3$ is not $\R$-factorizable.

In \cite{Sipacheva2023}, topological groups $\sg G4$ and $\sg H4$ are constructed with the same properties as the groups $\sg G3$ and $\sg H3$ and, additionally, the following condition is satisfied: the group $\sg G4\times \sg H4$ embeds closed in some strongly zero-dimensional group. The first and second arguments also work for the group $\sg G4\times \sg H4$, hence $\sg G4\times \sg H4$ is not $\R$-factorizable. In \cite{Sipacheva2023} a third argument was found why the group $\sg G4\times \sg H4$ is not $\R$-factorizable.
Thirdly, $\sg G4\times \sg H4$ is not strongly zero-dimensional and is embedded in some strongly zero-dimensional group, hence, by Theorem \ref{t:pog-defs:rdim}(1), $\sg G4\times \sg H4$ is not $\R$-factorizable.

%\section{Answers} \label{sec:answers} \input answers.tex

\bibliographystyle{elsarticle-num}
\bibliography{products_of_groups}
\end{document}